\documentclass[a4paper,10pt, twoside]{article}
\usepackage[english]{babel}
\usepackage{amsmath, amssymb,theorem}
\usepackage[all]{xy}
\usepackage{epsfig,color,rotating}
\xyoption{2cell}
\UseAllTwocells
\UseComputerModernTips
\usepackage[a4]{Martijn}
\newcommand{\note}{\textit}
\theoremheaderfont{\normalfont\bfseries}
\theorembodyfont{\normalfont\itshape}
\theoremstyle{change}
\newtheorem{Tm}{\normalfont\scshape Theorem}[section]
\newtheorem{Pp}[Tm]{\normalfont\scshape Proposition}
\newtheorem{Lm}[Tm]{\normalfont\scshape Lemma}
\newtheorem{Cr}[Tm]{\normalfont\scshape Corollary}

\theorembodyfont{\rmfamily}

\newtheorem{Df}[Tm]{\normalfont\scshape Definition}
\newtheorem{Ex}[Tm]{\normalfont\scshape Example}
\newtheorem{Rm}[Tm]{\normalfont\scshape Remark}

\newenvironment{Pf}{{\par\scshape Proof
}}{\hspace*{\fill}{\scshape QED}\\ \vspace*{0.5cm}\par}
\makeatletter
\renewcommand\subsection{\@startsection{subsection}{2}{\z@}%
                                     {-3.3ex\@plus -1ex \@minus -.2ex}%
                                     {2ex \@plus .2ex}%
                                     {\centering\normalfont\bfseries}}
\renewcommand{\@makecaption}[2]{%
\vspace{10pt}\hspace{.1\linewidth}
\parbox[l]{.8\linewidth}{\footnotesize{\scshape #1:} #2}%
\par
}
\makeatother

\newcommand{\RR}{\mathbb{R}}

\newcommand{\NN}{\mathbb{N}}

\newcommand{\QQ}{\mathbb{Q}}

\renewcommand{\phi}{\varphi}

\newcommand{\cobar}{\smash{\raise.5\baselineskip\hbox{\begin{turn}{180}$B$\end{turn}}}}
\newcommand{\pd}{\partial}

\newcommand{\End}{\mathrm{End}}
\newcommand{\Aut}{\mathrm{Aut}}
\newcommand{\id}{\mathrm{id}}

\renewcommand{\vert}{\textbf{v}}
\newcommand{\edge}{\textbf{e}}
\newcommand{\leg}{\textbf{l}}

\newcommand{\Opd}{\textsf{Opd}}

\newcommand{\Coopd}{\textsf{Coopd}}

\title{Coloured Koszul duality and \\ strongly homotopy operads}
\author{Pepijn van der Laan}
\begin{document}
\maketitle
\abstract{This paper proves Koszul duality for coloured operads and
  uses it to introduce strongly homotopy operads as a suitable
  homotopy invariant version of operads. It
  shows that $\QQ$-chains on configuration spaces of points in the
  unit disk form a strongly homotopy operad quasi isomorphic to the chains
  on the little disks operad.}

\section{Introduction}

Throughout this paper operads are operads in the category of dg vector
spaces over a field $k$ of characteristic 0.
\\~\\
In some situations the notion of operad is too restrictive. Think of
the following.
\begin{enumerate}
\item Given two quasi isomorphic operads $P$ and $Q$ there need not
  exist a quasi isomorphism $P\longrightarrow Q$ of operads.
\item Given an operad $P$, one usually can not
  transfer a strongly homotopy $P$-algebra structure from a dg vector
  space $W$ to a dg vector space $V$ using a map
  $\End_W\longrightarrow \End_V$, which from the operadic point of
  view would be the most natural thing to try.
\item The singular $k$-chains on configuration spaces of distinct
  ordered points in the unit disk in do not form an operad quasi
  isomorphic to the $k$-chains on the little disks operad in any
  straightforward manner, unless one uses Fulton-MacPherson
  compactification.
\end{enumerate}
The way in which this paper deals with these difficulties is by
defining a somewhat weaker version of operads, strongly homotopy
operds and morphisms between them. The
definition of a strongly homotopy operad is based on the analogy
between operads and associtive algebras advocated by
Ginzburg-Kapranov \cite{GinKap:Koszul}. In this analogy strongly
homotopy operads correspond to $A_\infty$-algebras (i.e. strongly
homotopy associative algebras). This paper shows that one can make the
analogy very precise
using Koszul duality for the $\NN$-coloured operad which as as
algebras non-symmetric pseudo operads. In fact one recovers the
associative algebra analogon when restricting to s.h. operads $P$ that as
collections  are concentrated in $P(1)$.

The main results can be summarized as follows.
\begin{enumerate}
\item Every quasi isomorphism of strongly homotopy operads admits a
  quasi inverse. Consequently, two augmented operads $P$ and $Q$ are quasi
  isomorphic iff there exists a quasi isomorphism $P\leadsto Q$ of
  strongly homotopy operads.
\item If $W$ and $V$ are two dg vector spaces and $i:V\longrightarrow
  W$,  $r:W \longrightarrow V$ and $H:W\longrightarrow W[1]$ are dg
  maps such that $H$ is a chain homotopy between $i\circ r$ and the
  identity on $W$, then there exists a morphism of strongly homotopy
  operads 
\[
\End_W \leadsto \End_V,
\]
 given by an explicit formula. This map is a quasi isomorphism if $r$
 and $i$ are quasi isomorphisms. If $P$ is an operad, and $W$ is a
 strongly homotopy $P$-algebra, this map can be used to transport the
 strongly homotopy $P$-algebra structure to $V$.
\item The $k$-chains on configuration spaces of ordered distinct
  points in the unit disk form a strongly homotopy operad quasi
  isomorphic to the $k$-chains on the little disks operad.
\end{enumerate}
Further applications of strongly homotopy operads related to
formal deformation theory for operads and their algebras, and
$L_\infty$-algebras can be found in my thesis \cite{Pep:Thesis}.

\subsection{Plan of the paper}

The preliminaries (Section \ref{Sec:Prelim}) fix some
notation. Section \ref{Sec:Coloured} briefly introduces coloured
operads, and then 
shows Koszul dality can be extended to coloured operads.
Section \ref{Sec:Nonsigmaopd} applies this to the $\NN$-coloured
operad $\mathrm{PsOpd}$ which has as algebras non-symmetric pseudo
operads, and gives an equivariant version of strongly homotopy
$\mathrm{PsOpd}$-algebras that defines strongly homotopy
operads. Finally, it conciders morphisms of strongly homotopy operads
and proves the first main result.
Section \ref{Sec:homotopyalgebras} proves the second main result, and considers
its implications for strongly homotopy $P$-algebras. The application
of these results to the operad $\mathrm{PsOpd}$ lead to the proof of
the third main result. 

\subsection{Acknowledgements}

I am grateful to Ieke Moerdijk and Martin Markl for enjoyable and
useful discussions. The research is part of my Ph.D. thesis
\cite{Pep:Thesis}, and was partly supported by Marie Curie 
Training Site Fellowships HPMT-CT-2001-00367 (Universit\'e Paris Nord
XIII) and HPMT-2000-00075 (Centre de Recerca Matem\`atica, Barcelona).

\section{Preliminaries}
\label{Sec:Prelim}
I work in the category of dg vector spaces over a field $k$ of
characteristic 0. If $V$ is a dg vector space, and $v\in V$ is an
homogeneous element, then its degree will be denoted by $|v|$. I use
the cohomological convention: the differential $d$ of the 
dg vector space $V$ is a map of degree $+1$. Let $V^n =
\{v\in V| |v|=n\}$ be the space of homogeneous elements of degree
$n$. Then $V[m]$ is the dg 
vector space with $(V[m])^n = V^{n-m}$. Later on I might
be a bit sloppy and leave out the `dg' since I only work with
differentially graded objects. Let $V$ and $W$ be (dg) vector
spaces. Recall that the symmetry $\tau$ of the tensor product involves
the natural signs $\tau:v\otimes w \longmapsto (-1)^{|v||w|}w\otimes
v$ on homogeneous elements.  I use the Koszul convention $f\otimes g
(x\otimes y) := (-1)^{|g||x|}f(x)\otimes g(y)$ for homogeneous maps
$f:V\longrightarrow V'$ and $g:W\longrightarrow W'$ applied to
homogeneous elements $x\in V$ and $y\in W$. In combination with the
shift $[m]$ this reduces the number of signs significantly.

By $S_n$ we denote the symmetric group on $n$ letters, and by 
$kS_n$ its group algebra which is the vector space spanned by the set
$S_n$ whose multiplication is the linear extension of multiplication
in $S_n$. If a group $G$ acts on a vector space $V$, the coinvariants of the
group action are denoted $V_{G}$ and the invariants by $V^{G}$.


\subsection{Operads}

A non-symmetric operad is a sequence
$\{P(n)\}_{n\geq 1}$ of (dg) vector spaces together with composition maps  
\[
\gamma:P(n)\otimes P(m_1)\otimes\ldots\otimes P(m_n) \longrightarrow
P(m_1+\ldots+ m_n),
\]
and an identity element $\id\in P(1)$. These structures satisfy the
usual associativity and identity axioms (cf. Getzler-Jones
\cite{GetzJon:Opd}, and Markl-Shnider-Stasheff \cite{MarShniSta:Opd}).

A collection $P$ is a sequence of vector spaces $\{P(n)\}_{n\geq 1}$
such that each $P(n)$ has a right $S_n$-module structure. An
(symmetric) operad is a collection $P$ together with an
non-symmetric operad structure on the sequence of vector spaces, and
composition is equivariant with respect to the $S_n$-actions in the
usual sense (cf. Getzler-Jones \cite{GetzJon:Opd}, and
Markl-Shnider-Stasheff \cite{MarShniSta:Opd}). Dually (in the sense of
inverting direction of arrows in the defining diagrams), one defines
(non-symmetric) cooperads. 

(Non-symmetric) pseudo operads, are the non-unital analogon of
(non-symmetric) operads. A non-symmetric pseudo operad $P$ is a
sequence of dg vector spaces $P$ together with dg maps
$\circ_i:P(n)\otimes P(m)\longrightarrow P(m+n-1)$ for $i=1,\ldots,
n$, which satisfy the appropriate associativity conditions. A
non-symmetric operad gives rise to a non-symmetric pseudo operad by
\begin{equation}\label{eq:circi}
p\circ_i q = \gamma(p;\id^{i-1},q,\id^{n-i}) 
\end{equation}
for $p\in P(n)$ and $q\in P(m)$. Pseudo operads are the equivariant
version of this, starting from a collection $P$. The category of
(non-symmetric) pseudo operads is equivalent to the category of
augmented operads. That is, operads $P$ such that the inclusion of the 
identity is split as a map of operads. Throughout this paper I assume
all operads except endomorphism operads to be augmented.
\\~\\
A graph $\eta$ consist of sets $\vert(\eta)$ of vertices, a set
$\edge(\eta)$ of internal edges, and a set $\leg(\eta)$ of external
edges or legs; together with a map that assigns to each edge a pair of
(not necessary distinct) vertices and a map that assigns to each leg a
vertex. To draw a graph, draw a dot for each vertex $v$, and for each
edge $e$ draw aline between the two vertices assigned to it, and for
each leg draw a line that in one end ends in the vertex assigned to
it. If $v\in \vert(\eta)$, denote by $\leg(v)\subset
\edge(v)\cup\leg(\eta)$ the set of legs and edges attached to
$v$ and call elements of $\leg(v)$ the legs of $v$. A morphism of
graphs consists of morphisms of vertices, edges, and legs compatible
with the structure maps. 

A connected graph $t$ is a tree if $|\vert(t)| = |\edge(t)|+1$. 
A rooted tree is a tree together with a basepoint $r\in\leg(t)$, the
root, and together with a bijection $\leg(t)\longrightarrow
\{0,\ldots,n_t\}$ that sends the root to 0, where $n_t =
|\leg(t)|-1$. In a rooted tree $t$, each of the sets $\leg(v)$ has a natural 
basepoint, the leg in the direction of the root. A planar tree is a
rooted tree together with for each $v\in\vert(t)$ a bijection
$\leg(t)\longrightarrow\{0,\ldots,n_v\}$ that sends the 
basepoint to 0, for $n_v := |\leg(v)|-1$. For any planar tree $t$ define
\[
C(t) := \bigotimes_{v\in\vert(t)}C(n_v).
\]
The free pseudo operad  $TC$  and the `cofree'
pseudo cooperad $T'C$ on a collection $C$ satisfy
\[
TC(n) =\underset{t}{\text{colim} } C(t), \qquad
T'C(n) = \lim_{t} C(t),
\]
where both limit and colimit are over the groupoid of planar trees
with $n$ external edges different from the root with isomorphisms of
rooted trees as maps. These maps need not preserve the planar
structure, but do preserve the labeling of the legs different from the
root in $\leg(t)$ by $1,\ldots, n$. The operad structure on $TC$ is
given by grafting trees, while the cooperad structure on $T'C$ is
given by cutting edges.

It is useful to be a bit more explicit on the arrows of the diagram
over which we take the (co)limit in defining $TC$ and $T'C$. Let
$\sigma:t\longmapsto t'$ be an isomorphism of rooted trees. For 
$v\in\vert(t)$ and $v'\in\vert(t')$, if $\sigma(v) = v'$ it
induces $C(n_v)\longrightarrow C(n_{v'})$. Define
\[
C(\sigma):C(t) = \bigotimes_{v\in t} C(n_v)\longrightarrow
\bigotimes_{v'\in t'} C(n_{v'}) = C(t'),
\]
as the tensor product over $v\in \vert(t)$ of these maps.
Note that $\sigma$ 
restricts to a bijection $\sigma|_{\leg(t)}:\leg(t)\longrightarrow
\leg(t')$, and that $\sigma$ being an isomorphism of rooted trees
implies compatibility of the labeling of the external edges of the
trees with $\sigma|_{\leg(t)}$.

If $P$ is an operad, then there is a natural differential $\pd_P$ on
$T'(P[-1])$. That is, $\pd_P$ is a square-zero coderivation of degree
+1. The resulting cooperad $BP = (T'(P[-1]),\pd)$ is the bar
construction on $P$. For more extensive background on this and on
(co)operads in general read Ginzburg-Kapranov \cite{GinKap:Koszul},
Getzler-Jones \cite{GetzJon:Opd}, and Markl-Shnider-Stasheff
\cite{MarShniSta:Opd}.

\section{Coloured Koszul duality}
\label{Sec:Coloured}
\subsection{Coloured operads} 

Denote by $\mathbf{n}$ the set $\{0,1,\ldots,n\}$ for $n \geq 0$, and
let $I$ be a set. An \note{$I$-coloured collection} (or \note{$I$-collection})
$P$ is a set
$\{P(n,i)\}_{(n,i:\mathbf{n}\rightarrow I)}$  of dg vector spaces
indexed by the sets $\mathbf{n}=\{0,1,\ldots,n\}$ for all $n\geq1$,
and by all surjections $\mathbf{n}\longrightarrow I$; together with a right
$S_n$-action on $\bigoplus_{i:\mathbf{n}\rightarrow I} P(n,i)$ such
that for $\sigma\in S_n$ the action satisfies $(P(n,i))\sigma\subset
P(n,i\sigma)$, where $i\sigma: \mathbf{n}\longrightarrow I$ is $i$
precomposed by the permution $\sigma$ applied to
$\{1,\ldots,n\}\subset \mathbf{n}$. The values of $i$ are called
labels. More particular, $i(0)$ is the output label, and
$i(1),\ldots,i(n)$ are the labels of the inputs $1,\ldots,n$. 

An \note{$I$-coloured pseudo operad} (or \note{$I$-pseudo operad}) is an $I$-collection
$P$, together with compositions
\[
\circ_l:P(n,i)\otimes P(m,j) \longrightarrow P(m+n-1, i\circ_lj) \qquad
\text{for } l\leq n, \text{s.t. } i(k)=j(0),
\]
( compare equation (\ref{eq:circi})) where
$i\circ_lj:\mathbf{(m+n-1)}\longrightarrow I$ satisfies 
\[
i\circ_lj(k) = \left\{\begin{array}{c @{\qquad\text{if }}c} i(k) &
  0\leq k<l\\
  j(k-l+1) & l\leq k <l+m \\ i(k-m) & k\geq l+m. \end{array}\right. 
\]
These data satisfy the compatibility relations for
$\circ_k$-operations of a pseudo operad (associativity, equivariance)
whenever these make sense.

An \note{$I$-coloured operad} (or \note{$I$-operad}) is an $I$-pseudo operad together
with for each $\alpha\in I$ an identity $\id_\alpha\in
P(\mathbf{1},\alpha)$, where $\alpha:\mathbf{1}\longrightarrow I$ is
the constant map with value $\alpha$. These identities act as units
with respect to any well defined composition.

Similarly, define $I$-pseudo cooperads, and $I$-cooperads by
inverting the arrows in the defining diagrams. 

\begin{Ex}
There is an obvious 1-1 correspondence between operads and
$*$-operads, where $*$ is the one-point set.

Let $I$ be a set. and Let $V = \{V_\alpha\}_{\alpha\in I}$ be a set of
vectorspaces. Denote by $Hom_k(-,-)$ the $k$-linear maps
(the internal Hom functor). Then $\End_V(n,i):= Hom_k(V_{i(1)}\otimes\ldots\otimes 
V_{i(n)}, V_{i(0)})$ defines a $I$-operad with respect to the
$S_n$-action on inputs and the obvious composition of maps where
$\phi \circ_l\psi$ uses the output of the map $\psi$ s the $l$-th
input of $\phi$. This $I$-operad is called the \note{endomorphism operad} of $V$. 

Let $P$ be an $I$-operad. A \note{$P$-algebra} $V$ is a set of vector spaces
$V = \{V_\alpha\}_{\alpha\in I}$ together with a morphism of
$I$-operads $P\longrightarrow \End_V$. 
\end{Ex}

Let $I$ be a set, and denote by $A_I$ the (non-unital) associative
algebra generated by genertors $[\alpha]$ for $\alpha\in I$ with the
multiplication $[\alpha] \cdot [\alpha'] = \delta_{\alpha\alpha'} \cdot[\alpha]$,
where $\delta$ is the Kronecker delta on the set $I$. For an
associtive algebra $A$, recall the definition of an $A$-pseudo operad
as a pseudo operad in the category of $A$-modules. The following is
now quite straightforward.
\begin{Ex}\label{Ex:AI}
Every $I$-coloured collection $P$ gives rise to a collection in the
collection of $A_I$-modules if we interpret
$\bigoplus_{i:\mathbf{n}\rightarrow I} P(n,i)$ as a decomposition in
eigenspaces of the left $A_I$ and right $A_I^{\otimes n}$-action with
eigenvalue 1. The left action of a generator $[\alpha]$ on $P(n,i)$ is
(again in terms of the Kronecker delta) $\delta_{\alpha\,i(0)}\cdot\id$ and
the right action of $[\alpha_1]\otimes \ldots \otimes [\alpha_n]$ on
$P(n,i)$ is $\delta_{\alpha_1\,i(1)} \cdot \ldots \cdot
\delta_{\alpha_n\,i(n)}\cdot\id$.
\end{Ex}
\begin{Pp}[Markl \cite{Mar:HomDiag}]\label{Pp:correspond}
There is a 1-1 correspondence between $I$-pseudo operads and
$A_I$-pseudo operads together with a decomposition 
\[
P(n) = \bigoplus_{i:\mathbf{n}\rightarrow I} P(n,i)
\]
that makes $P$ an $I$-collection (with the action of Example \ref{Ex:AI}
above) and such that for $p\in P(n,i)$ and $q\in P(m,j)$   
\[
p\circ_l q = 0 \qquad \text{if } j(0)\neq i(l).
\]
This correspondence describes $I$-operads as a full subcategory of
$A_I$-operads.
\end{Pp}

\subsection{Koszul duality for $I$-operads}

I assume the reader is familiar with Koszul duality for operads as
introduced in Ginzburg-Kapranov \cite{GinKap:Koszul} and its
description using cooperads in Getzler-Jones \cite{GetzJon:Opd}.
To prove that Koszul duality works for $I$-operads it suffices to
show that $I$-operads are closed under the relevant constructions in
the category of $A_I$-operads. 
\begin{Lm}\label{Lm:Ibarconst}
The bar construction $B_{A_I}$ from $A_I$-pseudo operads to
$A_I$-pseudo cooperads restricts to a functor $B_I$ from $I$-operds to
$I$-cooperads.  
\end{Lm}
\begin{Pf}
Let $P$ be an $I$-pseudo operad considered as an $A_I$-pseudo operad.
Recall that $B_{A_I}P(n)$ decomposes as a sum over trees with $n$
leaves with vertices labeled by elements of $P$. Each action of
$A_I$ corresponds to an input or output in $B_{A_I}P(n)$.It thus is
the action on the label of the vertex to which the corresponding 
leaf or root is attached. We get a decomposition of $B_{A_I}P(n)$ by
the generators $[\alpha]$ that do not vanish on these labels. 
\end{Pf}
An ideal $J$ of an $I$-operad $P$ is a sub $I$-collection $J$ of $P$ such that
$p\circ_lq\in J$ iff either $p$ or $q$ is an element of $J$. Denote
the free $I$-operad on a collection $E$ by $T_IE$.
An $I$-operad is called quadratic if it is of the from $T_IE/R$, where
$E(n,i) = 0$ if $n\neq 2$, and $R$ is an ideal generated by elements
in $\bigoplus_{i:\mathbf{3}\rightarrow I}T_IE(3,i)$. Quadratic
operads are naturally augmented.
\begin{Df}
The \note{Koszul dual $I$-cooperad} $P^\bot$ of an $I$-operad $P$ is
its Koszul dual as an $A_I$-operad. The Lemma below shows this is well
defined. An quadratic $I$-operad is Koszul if
$P^\bot\longrightarrow B_A(P)$ is a quasi isomorphism of
$A_I$-cooperads. The Koszul dual $I$-operad is $P^!=(P^\bot)^*\otimes
\Lambda$, the linear dual of $P^\bot$ tensored with the determinant
operad (cf. Getzler-Jones \cite{GetzJon:Opd}).
\end{Df}
\begin{Lm}\label{Lm:Ibot}
If $P$ is a quadratic $I$-operad, then the Koszul dual $P^{\bot}$ of
$P$ is an $I$-cooperad.
\end{Lm}
\begin{Pf}
Let $P=T_IE/R$ be a quadratic $I$-operad. The free $A_I$-cooperad
$T'(E[-1])$ is an $I$-cooperad by the same argument  on trees as
above. Moreover, by categorical generalities it is the free
$I$-coloured cooperad $T_I(E[-1])$ (under the correspondence of
Proposition \ref{Pp:correspond}).

The definition of $P^\bot$ as the kernel
of $T'(E[-1])\longrightarrow T'(R')$ where $R' =
T'(E[-1])(3)/R(3)[-2])$, assures that $P^\bot$ is an $I$-cooperad since
$R$ is an ideal.
\end{Pf}
Let $P$ be a quadratic operad.
The Koszul complex of a $P$-algebra $K$ is the cofree
$P^\bot$-coalgebra on the shifted vector space $K[-1]$, with the
natural differential obtained from the $P$-algebra structure on $K$ in
the sense of Ginzburg-Kapranov \cite{GinKap:Koszul}. It's homology is
denoted $H^P_*(K)$.
\begin{Tm}\label{Cr:colouredKoszul}
Let $P$ be a quadratic $I$-coloured operad.
\begin{enumerate}
\item The $I$-operad $P$ is Koszul iff $P^\bot\longrightarrow B_A(P)$
  is a quasi isomorphism of $I$-operads.
\item The $I$-operad $P$ is Koszul iff $P^!$ is Koszul
\item The homology $H^P_*(K)$ of the Koszul complex of a $P$-algebra
  $K$ vanishes for every free $P$-algebra $K$ iff $P$ is Koszul.
\end{enumerate}
\end{Tm}
\begin{Pf}
The result follows directly from the Lemmas \ref{Lm:Ibarconst} and
\ref{Lm:Ibot}, and Koszul duality 
for operads over a semi-simple algebra as proved in Ginzburg-Kapranov
\cite{GinKap:Koszul}.
\end{Pf}
\begin{Rm}
This article is devoted to one example of coloured Koszul
duality. More examples can be found in \cite{Pep:Thesis}. Coloured
Koszul duality is independently proved by 
Longoni and Tradler in preprint \cite{LonTra:Kos}.

Koszul duality has a nice interpretation in terms of the model category
of $I$-operads, the existence of which can be proved by the methods of
Berger-Moerdijk \cite{BerMoer:Model}. Namely, if $P$ is a Koszul
$I$-operad, then $\cobar_I(P^{\bot})\longrightarrow P$ gives a concise
cofibrant replacement for augmented operads $P$ in this model
category, where $\cobar_I$ is the cobar construction from
$I$-cooperads to $I$-operads (the dual construction to $B_I$ in Lemma
\ref{Lm:Ibarconst}).
\end{Rm}

\section{Strongly homotopy operads}\label{Sec:Nonsigmaopd}

\subsection{An operad of non-symmetric pseudo operads}

%
\begin{Df}
Define an $\NN$-operad $\mathrm{PsOpd}$ as follows.
As an $\NN$-collection, $\mathrm{PsOpd}(n,i)$ is spanned by planar
rooted trees $t$ with $n$ vertices numbered $1$ up to $n$, that
satisfy $|\leg_t(k)|-1 = i(k)$ for $k=1,\ldots,n$, and $i(0) = |\leg(t)|-1$.
Composition $s\circ_k t$ is defined by
replacing vertex $k$ in $s$ by the planar rooted tree $t$ (cf. Figure
\ref{Fig:circPsOpd}). More precisely, $s\circ_k t$ has vertices
$\vert(s)-\{k\} \cup \vert(t)$, and edges $\edge(s)\cup \edge(t)$,
where the elements of $\leg_s(k)$ necessary to define the edges of $s$
are interpreted as elements of $\leg(t)$. This is well defined since
the planar structure gives a natural isomorphism between $\leg(t)$ and
$\leg_s(k)$.   
\begin{figure}[!ht]
\begin{center}
\input{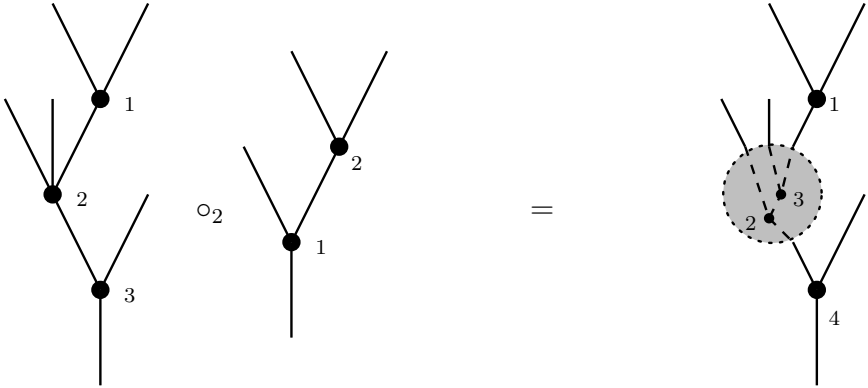}
\caption{Composition $\circ_2$ in $\mathrm{PsOpd}$: vertex
$2$ of the left tree is replaced by a tree with matching number of
legs.\label{Fig:circPsOpd}} 
\end{center}
\end{figure}
\end{Df}
\begin{Pp}\label{Tm:PsOpdonsigtree}
The $\NN$-operad $\mathrm{PsOpd}$ is a quadratic
$\NN$-operad. Algebras for $\mathrm{PsOpd}$ are
non-symmetric pseudo operads.
\end{Pp}
\begin{Pf}
Every planar rooted tree can be constructed from
2-vertex trees by compositions in $\mathrm{PsOpd}$, adding
one edge at a time. Denote
by $(\mathbf{m}\circ_i\mathbf{n})$ the 2-vertex planar rooted tree
with the root vertex having legs $\{0,\ldots,m\}$, and the other
vertex having legs $\{0,\ldots,n\}$. The unique internal edge connects
leg $i$ of the root vertex to leg $0$ of the other vertex. These generators
satisfy the quadratic relations
\begin{equation}\label{Eq:nonsigmapsopdrel}
(\mathbf{k}\circ_j (\mathbf{m}\circ_i \mathbf{n})) = \left\{
\begin{array}{r @{\qquad} l }
((\mathbf{k}\circ_{i}\mathbf{n}) \circ_{j+n-1} \mathbf{m}) &  \text{if }i<j \\
(\mathbf{k} \circ_{j} (\mathbf{m}\circ_{i-j+1} \mathbf{n})) & \text{if }j\leq i
< j+ m \\
((\mathbf{k} \circ_{i-m-1} \mathbf{n})\circ_{j} \mathbf{m})) &  \text{if }j\leq i+m
\end{array}\right.
\end{equation}
These generators and relations define a quadratic
$\NN$-operad with free $S_n$-actions and non-symmetric pseudo operads
as algebras, as follows from the definition. Denote this quadratic operad $TE/R$. 

To identify the two 1-reduced operads
$\mathrm{PsOpd}$ and $TE/R$ it suffices to identify the
free algebras on 1 generator in each colour since both $\NN$-operads
have free $S_n$-actions. Recall that the free 
non-symmetric pseudo operad on $A= \{A_n\}_{n\in\NN}$ is given as
$\bigoplus_t A(t)$, where the sum is over planar trees (cf. Loday
\cite{Loday:dialgebras}, Appendix B). Hence the free algebras are
isomorphic.
\end{Pf}

\subsection{Koszul duality for $\mathrm{PsOpd}$}

\begin{Tm}\label{Tm:PsOpdKoszul}
The $\NN$-coloured operad $\mathrm{PsOpd}$ of non-symmetric pseudo operads is a
self dual Koszul $\NN$-coloured operad. 
\end{Tm}
\begin{Pf}
Write $\mathrm{PsOpd} = TE/R$ as in the proof of the previous
result. We compute the Koszul dual operad $\mathrm{PsOpd}^! =
T(E^*)/R^\bot$,
where $R^\bot$ is the ortogonal complement of $R$ with respect to the
pairing of $TE^*$ and $TE$ defined as the extension of the piring of
$E^*$ and $E$ twisted by a sign (cf. Ginzburg-Kapranov
\cite{GinKap:Koszul}).

The dimension of $R(3)$ is exactly half the dimension
of $\mathrm{PsOpd}(3)$, since the associativity relations
divide the basis elements of $\mathrm{PsOpd}(3)$ in
pairs which satisfy a non-trivial relation. Observe that the dual
relations $R^{\bot}(3)$ certainly are contained in  the ideal
generated by 
\[
(\mathbf{k}\circ_j (\mathbf{m}\circ_i \mathbf{n})) = \left\{
\begin{array}{r @{\qquad} l }
((\mathbf{k}\circ_{i}\mathbf{n}) \circ_{j+n-1} \mathbf{m}) &\text{if } i<j \\
-(\mathbf{k} \circ_{j} (\mathbf{m}\circ_{i-j+1} \mathbf{n}))&
\text{if }j\leq i < j+ m \\
((\mathbf{k} \circ_{i-m-1} \mathbf{n})\circ_{j} \mathbf{m}))&\text{if }j\leq i+m
\end{array}\right.
\]
By a dimension argument these relation must exactly all the
relations. Then a base change shows that $(\mathrm{PsOpd})^!$ is
isomorphic to  $\mathrm{PsOpd}$. The base change is given by
multiplying a basis element corresponding to a planar rooted tree $t$
with the sign $(-1)^{c(t)}$, where $c(t)$ is the number of internal
\textit{axils} of $t$. That is, the number of distinct subsets
$\{v,w,u\}\subset\vert(t)$ such that two of the three vertices are
direct predecessors of the third. This shows that
$\mathrm{PsOpd}$ is self dual. 
 
Let $P$ be a non-symmetric pseudo operad. The $\mathrm{PsOpd}$-algebra 
homology complex of $P$ is as a sequence of graded vector spaces the
free non-symmetric pseudo cooperad on $P[-1]$,
\[
C_*^{\mathrm{PsOpd}}(P) = F'_{\mathrm{PsOpd}^\bot}(P)
= \bigoplus_{t \mathrm{\ planar}}
\bigotimes_{v\in t} P(\leg_t(v))[-1].
\]
The differential is given by contracting edges using the
$\circ_i$-compositions in $P$. In other words, this complex
is the non-symmetric bar construction $B_{\not\Sigma}P$ (cf. Loday
\cite{Loday:dialgebras}, appendix B). 
The Theorem follows since the homology of this
complex vanishes in the case where $P=T_{\not\Sigma}C$, the free
non-symmetric operad on $C$.
\end{Pf}
\begin{Rm}
Theorem \ref{Tm:PsOpdKoszul} invites the reader to a conceptual
excursion. As explained in the proof, the homology complex 
\[
C_*^{\mathrm{PsOpd}}(P) =
(F'_{\mathrm{PsOpd}^\bot}(P),\pd)
\]
of a non-symmetric
pseudo operad $P$ is the non-symmetric bar complex of $P$.
This shows how bar/cobar duality for non-symmetric
operads  is an example of Koszul duality for the coloured operad 
$\mathrm{PsOpd}$. The non-symmetric bar construction $B_{\not\Sigma}P$ of a
non-symmetric operad is nothing but the $\mathrm{PsOpd}$-algebra
complex of $P$, computing the $\mathrm{PsOpd}$-algebra
homology of the algebra $P$. 
\end{Rm}
\label{Df:htpynonsymPsOpd}
Let $P = \{P(n)\}_{n\geq 0}$ be a sequence of vector spaces. 
The formalism of Koszul duality defines a \note{strongly homotopy
  $\mathrm{PsOpd}$-algebra} (or a
\note{s.h. $\mathrm{PsOpd}$-algebra}) is a sequence of vector spaces
$P$, together with a square zero coderivation $\pd$ of the `cofree'
$\mathrm{PsOpd}^\bot$-coalgebra on $P$ of cohomological degree +1
(compare Ginzburg-Kapranov \cite{GinKap:Koszul}).

\begin{figure}[ht!]
\centering
\input{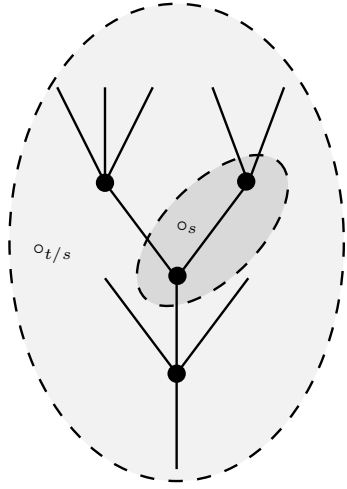}
\caption{One summand of Equation (\ref{eq:pdsquare}): $\circ_s$
contracts the darker part, covering subtree $s$ of $t$ and
$\circ_{t/s}$ contracts the remaining tree.\label{Fig:squarezero}} 
\end{figure}

For a planar rooted tree $t$, recall 
$P[-1](t) = \bigotimes_{v\in\vert(t)}P_{n_v}[-1]$, where $n_v=|\leg(v)|-1$.
A strongly homotopy $\mathrm{PsOpd}$-algebra structure on $P$ is
determined by operations 
\[
\circ_t:P[-1](t)\longrightarrow P[-1],
\] 
one for every planar rooted tree $t$.
The condition on $\pd^2 = 0$ on the differential is equivalent to a
sequence of relations on these 
operations. For each planar rooted tree $t$, we obtain a relation of the form
\begin{equation}\label{eq:pdsquare}
\sum_{s\subset t}\pm(\circ_{t/s})\circ(\circ_{s}) = 0,
\end{equation}
where the sum is over (connected) planar subtrees $s$ of $t$ and $t/s$ is the
tree obtained  from $t$ by contracting the subtree $s$ to a point, and
the signs involved are induced by a choice of ordering on the vertices
of the planar rooted trees $t$ and $s$ in combination with the Koszul
convention. Here a connected planar subtree is a subset of vertices
together with all their legs and edges such that the graph they
constitute is connected. One term of the sum is illustrated in 
Figure \ref{Fig:squarezero}.

%
%

%
%

\subsection{Strongly homotopy operads}

The s.h. $\mathrm{PsOpd}$-algebras described above are not quite what
we need, since these do not consider the symmetric group actions on
collections. 

\begin{Df}\label{Df:equivarianthtpyalg}
Let $P = \{P(n)\}_{n\in \NN}$ be a collection such that the
vector spaces $P(n)$ form a s.h. $\mathrm{PsOpd}$-algebra.
Let $t$ and $t'$ be planar rooted trees. If $\sigma:t\longmapsto t'$
is an isomorphism of the underlying rooted trees, then $\sigma$
induces $\sigma:P(t) \longmapsto P(t')$ through
the maps of $\Aut({\leg(v)})$-modules in the tensor factors
of $A(t)$, and it induced $\leg(\sigma): \leg(t)
\longrightarrow \leg(t')$ and consequently a map of
$\Aut(\leg(t))$-modules $\leg(\sigma):P(\leg(t)) \longrightarrow
P(\leg(t'))$. (Recall that the planar structure of $t$
induces a natural identification of $P(\leg(t))$ with $P(n)$, where $n
= |\leg(t)|-1$.)

Call a s.h. $\mathrm{PsOpd}$-algebra $P$
\note{equivariant}, if for every planar rooted tree $t$, and every
automorphism $\sigma$ as above, satisfies
\[
\leg(\sigma)\circ(\circ_t) = (\circ_{t'})\circ\sigma.
\]  
A \note{strongly homotopy operad} (or \note{s.h. operad}) is an
equivariant strongly homotopy $\mathrm{PsOpd}$-algebra.
\end{Df}
\begin{Rm}
Recall that $T'(P[-1])(n) = \lim_{t} P(t)$, where the
limit is over the groupoid of planar rooted trees
with $n$ leafs different from the root. A differential $\pd$ on
$F'_{\mathrm{PsOpd}^{\bot}}(P)$ defined by maps
$\circ_t:P(t)\longrightarrow P(\leg(t))$ induces maps 
on the limit $\lim_t P(t)\longrightarrow P(\leg(t))$ iff $\pd$ is
equivariant (i.e. defines a strongly homotopy operad). In that case it
defines a differential on $T'(P[-1])$. I use
notation
\[
BP = (T'(P[-1]),\pd),
\]
where $\pd$ denotes the induced differential.
\end{Rm}
\begin{Ex}
The  bar construction makes operads a
special case of operads up to homotopy, as is suggested by the
notations $BP$. The trees $t$  with $|\vert(t)| =
1$ define the internal differential, and the trees with $|\vert(t)| =
2$ the compositions $\circ_i$. The $\circ_t$ operations vanish if
$|\vert(t)|\geq 3$. The conditions on the $\circ_t$-operations
translate into the operad axioms. Operads are exactly
s.h. operads such that $\circ_s$ vanishes if $\vert(s) \geq 3$. 
\end{Ex}
\begin{Ex}
Interpret operads up to homotopy as a generalisation of operads where
one needs `higher homotopies' that measure the failure of
associativity of the $\circ_i$ operations. 
I dwell a bit on this interpretation:
Let $P$ be a s.h. operad.  When $|\vert(t)| = 1$, $\circ_t$ defines an
internal differential on $P(\leg(t))$. When $|\vert(t)| = 2$, the
operation $\circ_t$ defines a circle-$i$ operation as in the
definition of an (pseudo) operad. In general these operations
need no longer be associative. If $\circ_s$ does not vanish for 
$|\vert(s)| = 3$, then Equation (\ref{eq:pdsquare}) expresses (if
|$\vert(s)|=3$) that
$\circ_s$ serves an a homotopy for associativity as follows. Denote
the internal differential by $d$ and the two contractions of the
internal edges $e$ or $e'$ of $s$  by $\circ_e$ and $\circ_{e'}$ that
correspond to operad compositions. The formula
\[
(\circ_e)\circ(\circ_{e'}) - (\circ_{e'})\circ(\circ_e) = d\circ
(\circ_s) + (\circ_s)\circ d
\]
shows that associativity of the $\circ_e$ compositions holds up
to the homotopies $\circ_s$ with $|\vert(s)| = 3$.
More explicitly (with the signs), for a linear tree labelled with elements
$p,q,r$ in $P$ we have
\[
\begin{split}
(&p\circ_e q)\circ_{e'} r - p\circ_e(q\circ_{e'} r) \\
&= d(\circ_s(p,q,r))
+ \circ_s(dp,q,r) + (-1)^{|p|}\circ_s(p,dq,r) +
(-1)^{|p|+|q|}\circ_s(p,q,dr).  
\end{split}
\]
Consequently, if $P$ is an s.h. operad, then
the cohomology $H^*P$ with respect to the internal differential $d$
is a graded operad. The $\circ_i$-compositions are induced by the
operations $\circ_t$ for trees $t$ with 2 vertices.
\end{Ex}

\subsection{Homotopy homomorphisms}

\label{Sec:HomotopyHomomorphisms}
In the spirit of Koszul duality, we define a
\note{homotopy homomorphism of homotopy
$\mathrm{PsOpd}^\bot$-algebras} to be a morphism 
of cofree $\mathrm{PsOpd}^\bot$-coalgebras compatible with the
differentials.
Such a morphism is a quasi isomorphism 
if the underlying map of vector spaces is a quasi isomorphism.  
Recall that by the moves of Markl \cite{Mar:HomHom}
(we need to extend the theory to coloured operads but this is no
problem) such a quasi isomorphism has a quasi inverse.
A morphism of homotopy $\mathrm{PsOpd}$-algebras
$\phi:A\leadsto B$ is completely determined by its restrictions 
\[
\phi_t:(A[-1])(t) \longrightarrow B(\leg(t))[-1].
\]
The condition that 
$\phi$ is compatible with the differential can be described in terms of
conditions about compatibility with the $\circ_t$ operations:
\begin{equation}\label{Eq:homotopymorphism}
\sum_{s\subset t} \pm \phi_{t/(s)} \circ
(\circ_{s}) =  \sum_{n, s_1,\ldots,s_n\subset t} \pm 
(\circ_{t/(s_1,\ldots,s_n)})\circ(\phi_{s_1}\otimes\ldots\otimes \phi_{s_n}).
\end{equation}
where the sum in the left hand side is over subtrees, and the sum in
the right hand side for each $n$ is over $n$-tuples of (connected) subtrees of $t$ with
disjoint sets of internal edges that together cover all vertices of
$t$. The $\pm$ is the sign is induced by the Koszul 
convention.
\begin{Df}
A \note{homotopy homomorphism} $\phi:A\leadsto B$ of equivariant homotopy
$\mathrm{PsOpd}$-algebras is \note{equivariant} if for any
planar rooted trees $t$ 
and $t'$ and any isomorphism $\sigma:t\longrightarrow t'$ of the
underlying rooted trees, the equation  
\[
\leg(\sigma)\circ(\phi_t) = \phi_{t'}\circ\sigma
\]
is satisfied.
A \note{morphism of operads up to homotopy} $\phi:P\leadsto Q$  is an
equivariant homotopy homomorphism $\phi:P\leadsto Q$ of equivariant
homotopy $\mathrm{PsOpd}$-algebras. An equivariant s.h. morphism
$\phi:P\leadsto Q$ induces a morphism of dg cooperads
$\phi:BP\longrightarrow BQ$.


Note that $\phi$ is determined by
maps $\phi_t:(P[-1])(t)\longrightarrow P[-1]$.   
A \note{homotopy quasi isomorphism} is a homotopy homomorphism such
that the morphism $\phi_\bullet$ of dg collections is a
isomorphism in cohomology. Here $\phi_\bullet$ stands for the
restriction of $\phi$ to the 1-vertex trees.
\end{Df}

\begin{Tm}\label{Tm:quasiinv}
Let $P$ and $Q$ be s.h. operads, and let $\phi:P\leadsto Q$ be a
quasi isomorphism of s.h. operads. Then there exists a quasi inverse
$\psi:Q\leadsto P$ to $\phi$.
\end{Tm}
\begin{Pf}
A homotopy quasi
isomorphism of operads $\phi$ has a quasi 
inverse as a morphism of homotopy
$\mathrm{PsOpd}$-algebras 
(cf. Markl \cite{Mar:HomHom}). Let $\psi$
denote this quasi inverse. This quasi inverse can be symmetrised as
follows. Let $t$ be a planar rooted tree with $n$ vertices. Define
\[
\psi'(t) = \frac{1}{|\Aut(t)|}
\sum_{\sigma\in\Aut(t)}\psi_{\sigma(t)}\circ \sigma. 
\]
Since $\phi$ is equivariant, $\psi'$ still is a quasi inverse to
$\phi$. Moreover, for $\tau\in \Aut(t)$
\[|\Aut(t)| \cdot\psi_{t}') \circ \tau =
\sum_{\sigma\in\Aut(t)}\psi_{\sigma(t)}\circ \sigma \circ
\tau
= \sum_{\sigma'\in\Aut(t)}\psi_{\sigma'\circ\tau(t)}\circ \sigma',
\]
where we use $\sigma' = \sigma\circ\tau^{-1}$ to compare the sums.
Then $\psi'$ is an equivariant quasi inverse to $\phi$. 
\end{Pf}

\begin{Cr}\label{Cr:qiresult}
Two augmented operads $P$ and $Q$ are quasi isomorphic iff there
exists a quasi isomorphism $P\leadsto Q$ of operads up to homotopy. 
\end{Cr}
\begin{Pf}
By definition $P$ and $Q$ are quasi isomorphic iff there exists a
sequence of quasi isomorphisms of augmented operads
$P\longleftarrow \cdots \longrightarrow Q$. The previous theorem can
be applied to make all arrows point in the same direction if we allow
s.h. maps.

On the other hand, if there exists an s.h. quasi isomorphism
$P\leadsto Q$, then the bar-cobar adjunction (cf. Getzler-Jones
\cite{GetzJon:Opd}) gives a strict quasi isomorphim $\cobar(B
P)\longrightarrow Q$, where $\cobar(C)$ denotes the cobar construction
on a cooperad $C$. Moreover, there exists a natural quasi isomorphism
$\cobar(BP)\longrightarrow P$.
\end{Pf}
%
%
\section{Homotopy Algebras}
\label{Sec:homotopyalgebras}
%
%

\subsection{Endomorphism operads}

This section constructs  homotopy homomorphisms between
endomorphism operads, some even compatible with the identity. Well known
boundary conditions  turn up naturally in this context
(compare Huebschmann-Kadeishvili \cite{HuebKad:Model}). 

\begin{Df}
An s.h. operad is \note{strictly unital} if there exists an element
$\mathrm{id} \in P(1)$ that is a left and right identity with respect to
the $\circ_t$ operations where $|\vert(t)| = 2$ and such that the
other compositions $\circ_t$ vanish when applied to $\mathrm{id}$ in one
coordinate.
A homotopy homomorphism $\phi$ of two strictly unital operads up to homotopy
is \note{strictly unital} if the underlying
morphism $\phi_{\bullet}$ of collections preserves the identity, and if for
$|\vert(t)|>1$, the map $\phi(t)$ vanishes when applied to $\mathrm{id}$
in one coordinate. 
\end{Df}
Let $V$ and $W$ be dg vector spaces. $V$ is a \note{strict deformation retract} of
$W$ if there exist an inclusion $i:V\longrightarrow W$ and a retraction
$r:W\longrightarrow V$ such that both $i$ and $r$ are dg maps, $r\circ i
= \mathrm{id}_V$, and there exists a chain homotopy $H$  between $i\circ
r$ and $\id_W$, satisfying the boundary conditions $H\circ i = 0$,
$r\circ H = 0$, and $H\circ H = 0$.
\begin{Tm}\label{Lm:Endhomhom}
Let $V$ and $W$ be dg vector spaces. Let $i:V\longrightarrow
W$ and $r:W\longrightarrow V$ be dg linear maps, and $H:W
\longrightarrow W[1]$ a chain homotopy between $i\circ r$ and $\id_W$.
\begin{enumerate}
\item
There exists a (non-unital) homotopy homomorphism $\phi:\End_W\leadsto
\End_V$ (defined by the Formula (\ref{eq:Endhomohomform}) below). 
\item
If $i$ and $r$ are quasi isomorphisms,
then $\phi$ is a quasi isomorphism. 
\item
If the data above make $V$ a strict deformation retract of $W$, then
$\phi$ is strictly unital.
\end{enumerate}
\end{Tm}
\begin{Pf}
The map $\phi_\bullet$ corresponding to 1 vertex trees is
$f\longmapsto r\circ f \circ i^{\otimes n}$ for $f\in \End_V(n)$. This
proves the second part of the Theorem.
Define an alternative composition $\hat \gamma$ on $\End_V$ by $f\hat\circ_i
g = f \circ_1 H' \circ_i g$, where $H'(x) = (-1)^{|x|}H(x)$. This
composition makes $\End_V$ a pseudo operad. For a planar rooted tree $t$, the
map 
\begin{equation}
\label{eq:Endhomohomform}
\phi(t) = \phi_\bullet\circ \hat \gamma_t,
\end{equation}
where $\hat\gamma_t: \End_W(t) \longrightarrow \End_W(\leg(t))$ is
the composition based on $\hat\gamma$. This is visualised in Figure
\ref{Fig:nieuwecomp}.  
\begin{figure}[ht!]
\centering
\input{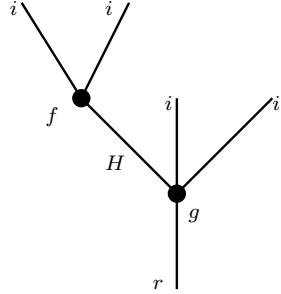}
\caption{The map $\phi(t)(f,g,h) = r\circ g \circ ((H' \circ f
\circ i^{\otimes 2}),i^{\otimes 2})$, represented by a tree with
labelled internal and external edges.\label{Fig:nieuwecomp}} 
\end{figure}

It remains to check Formula (\ref{Eq:homotopymorphism}).  
For a fixed tree $t$ this reduces to
\begin{equation}\label{Eq:Endhomohomo}
\sum_{e\in\edge(t)} (\circ_e)\circ (\phi(t^e)\otimes \phi(t_e)) +
d\circ\phi(t) = \sum_{e\in\edge(t)}\phi(t/e) \circ (\circ_e)
+\phi(t)\circ d. 
\end{equation}
The argument that this hold is the following. 
Since $r$ and $i$ commute with the differential $d$, and the internal
differentials act as a derivation with respect to composition of
multi-linear maps, the formula follows from the equalities $d\circ H +
H\circ d = \id - i \circ r$ applied to the summand for each edge
$e$. This shows part \textit{(i)}.

Assume the conditions of \textit{(iii)}. To assure that $\phi_\bullet$
preserves the identity, use $r\circ i = 
\id_V$. The conditions on compositions with $H$ assure that higher
operations applied to the identity vanish.
\end{Pf}
\begin{Rm}
Let us take a closer look at the proof above. Since the
cancellation of terms is local with respect to the geometry of the tree
$t$ (i.e. cancellation per edge), it suffices to check the signs for a
tree with one edge as in Figure \ref{Fig:nieuwecomp}. Let us do the
calculation with the signs for this tree. We leave out the
pre-composition with $i$ and post composition with $r$ in the final
terms. The usual degree of $f$ is denoted by $|f|$. 
The left hand side of Equation (\ref{Eq:Endhomohomo}) reads
\[
g\circ_k i \circ r \circ f + d\circ g\circ_k H\circ f +
(-1)^{|f|+|g|+1}g\circ_k H\circ f \circ d.
\]
The right hand side equals
\[
\begin{split}
&g\circ_k f  + (-1)^{|g|} g\circ d \circ_k
H \circ f + (-1)^{|g|+1} g\circ_k H\circ d\circ f \\
&+  d\circ g\circ_k H \circ f + (-1)^{|f|+|g|+1} g\circ_k H\circ f\circ d.
\end{split}
\]
To obtain the signs, note that we have a sign from moving $d$ in,
and note that these signs are with respect to the shifted grading on
$\End_V$ and $\End_W$, while the sign  in $d(f) = d\circ f + (-1)^{|f|}
f\circ d$ is with respect to the usual grading. The signs are correct
if we replace $H$ by $H'(x) = (-1)^{|x|} H(x)$.
\end{Rm}


\subsection{Homotopy $Q$-algebras}

I already discussed homotopy algebras for Koszul operads. This
section discusses the more general approach to homotopy algebras. It
shows how operads up to homotopy can be used to give a different
interpretation of the usual definition.

\begin{Df}\label{Df:Htpyalgdef}
Let $Q$ be an augmented operad. A \note{homotopy $Q$-algebra}
structure on a dg 
vector space $V$ is a homotopy homomorphism
$Q\leadsto\End_V$. 

Recall that this induces a map of cooperads $BQ\longrightarrow
B\End_V$. To such a morphism corresponds 
by the bar/cobar adjunction a morphism of operads $\cobar (BQ)\longrightarrow
\End_V$, where $\cobar:\Coopd\longrightarrow \Opd$ is the cobar
construction as in the proof of Corollary \ref{Cr:qiresult}. Moreover,
$\cobar B(Q)$ is a cofibrant replacement of $Q$ in the model category
of operads (cf. the proof of Corollary \ref{Cr:qiresult} for the
notation). This explains the terminology.
\end{Df}
\begin{Pp}\label{Tm:Linftymorphofdeformcplx}
Let $Q$ be an augmented operad.
\begin{enumerate}
\item
let $W$ be a homotopy $Q$-algebra, and $W$ a dg vector space. If
$i:V\longrightarrow W$, $r: W\longrightarrow V$ are quasi
isomorphisms, and $H:r\circ i \sim \id_W$, then 
$V$ has the structure of a homotopy $Q$-algebra such that the induced
maps in cohomology $H(r)$ and $H(i)$ are isomorphisms of $Q$-algebras. 
\item
Let $V$ be a homotopy $Q$-algebra, and let $W$ be a dg vector space. If
$i:V\longrightarrow W$, $r: W\longrightarrow V$ are quasi
isomorphisms, and $H:r\circ i \sim \id_W$, then 
$W$ has the structure of a homotopy $Q$-algebra such that $H(r)$ and
$H(i)$ are isomorphisms of $Q$-algebras. 
\end{enumerate}
\end{Pp}
\begin{Pf}
Suppose that $W$ is a homotopy $Q$-algebra. 
Recall that we constructed from the data in the Theorem a quasi isomorphism
$\End_W\leadsto \End_V$ in Theorem \ref{Lm:Endhomhom}. The composition 
\[
BQ\longrightarrow B\End_W \longrightarrow B\End_V
\]
defines the desired homotopy homomorphism $Q\leadsto \End_V$, where
the map $BQ\longrightarrow B\End_W$ is the map defined by the 
homotopy $Q$-algebra structure on $W$, which proves \textit{(i)}.

Suppose that $V$ is a homotopy $Q$-algebra. 
The quasi isomorphism $\End_W\leadsto \End_V$ has a quasi
inverse (by Theorem \ref{Tm:quasiinv}), and thus we can construct the
composition $BQ\longrightarrow B\End_V \longrightarrow B\End_W$, which
defines a homotopy $Q$-algebra structure on $W$.  
\end{Pf}

\begin{Rm}\label{Rm:ColCase}
Observe that all the results can be generalised to coloured operads: 
An \note{strongly homotopy $I$-operad} $P$ is an $I$-collection $P$ together
with a differential $\pd$ on the `cofree' pseudo $I$-cooperad
$T'_I(P[-1])$. If we denote $B_IP = (T_I'(P[-1]),\pd)$, we can define a
homotopy homomorphism $P\leadsto Q$ of s.h. $I$-operads as a 
morphism of $I$-cooperads $B_I P\longrightarrow B_I
Q$. Notably, for
sequences of vector spaces $V=\{V_\alpha\}_{\alpha\in I}$ and
$W=\{W_\alpha\}_{\alpha\in I}$ such that for each $W_\alpha$ and
$V_\alpha$ we have $i_\alpha$, $r_\alpha$ and $H_\alpha$ as in the
second part of the Proposition above, we can find a quasi isomorphism
$\End_W\leadsto \End_V$, which yields the analogue of Proposition
\ref{Tm:Linftymorphofdeformcplx} for algebras over $I$-operads. 
\end{Rm}

\subsection{Example: configuration spaces}

Let $D_2$ be the \note{operad of little disks}. That is, $D_2$ is the 
topological operad such that $D_2(n)$ is the space of ordered $n$-tuples 
of disjoint embedding of the unit disk $D_2$ in $D_2$ that preserve 
horizontal and vertical directions. The operations $\circ_k$ are defined 
by compositions of embeddings.  

Let $F(n)$ denote the \note{configuration space} of $n$ distinct
ordered points in the open unit disk in $\RR^2$. Thus $F(n)$ is the
$n$-fold product of the unit disk with the (sub)diagonals cut
out. Consider $F = \{F(n)\}_{n\geq 1}$ as a collection with respect
to permutation of the order of the points.

For a topological space $X$, denote by $S_*(X)$ the singular 
$k$-chain complex on $X$ with coefficients in $k$.

\begin{Tm}\label{Tm:confchain}
The singular $\QQ$-chains $S_*(F)$ on configuration spaces form an operad 
up to homotopy quasi isomorphic (in the sense of Proposition
\ref{Tm:Linftymorphofdeformcplx}) to the operad $S_*(D_2)$ of singular
$\QQ$-chains on the little disks operad.  
\end{Tm}
\begin{Pf}
We first sketch the line of argument.
We construct an $S_n$-equivariant homotopy between the little disks
and the configuration spaces. It then follows that $S_*(F)$ is a
homotopy algebra for $\mathrm{PsOpd}$ homotopy equivalent
to $S_*(D_2)$. Since the homotopy algebra $S_*(F)$ is equivariant,
$S_*(F)$ is an s.h. operad. This is based on the
observation in Remark \ref{Rm:ColCase} that the all results go through
for coloured operads.

Then there exists an inclusion $i:F(n)\longrightarrow 
D_2(n)$ and a retraction $r:D_2(n)\longrightarrow F(n)$ such that 
$\id\sim i\circ r$ by a homotopy $H$, and $r\circ i = \id$.
Consider points in $D_2(n)$ as given by a $n$-tuple $(x_1,\ldots,x_n)$ 
of points in the interior of $D_2$ and an $n$-tuple $(r_1,\ldots,r_n)$ 
of radii, and a point in $F(n)$ by a $n$-tuple $(x_1,\ldots,x_n)$ of 
points in the interior of $D_2$.  One might take the retraction $r$  by 
defining all radii in $r(x_1,\ldots,x_n)$ equal to 
\[
\frac{1}{3}(\mathrm{min}(\{|x_i-x_j| \quad (i\neq j)\}\cup \{ 1-|x_i|\})).
\]
(The map $r$ is not smooth but only continuous.)
A homotopy $H$ between $i\circ r$ and the identity is readily defined
by drawing a tube of configurations with the two configurations at the
boundary disks, connection the little disks by straight lines.
(cf. Figure \ref{Fig:tube}).
\begin{figure}[!ht]
\begin{center}
\input{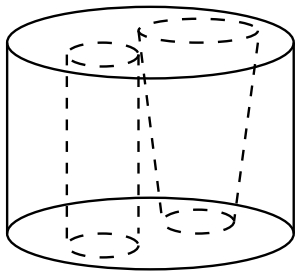}
\caption{Construction of the homotopy $H$. The tubes do not intersect
since the centres of the disks are fixed. \label{Fig:tube}}
\end{center}
\end{figure}

The homotopy $H$ induces a chain 
homotopy between $S_*(i)\circ S_*(r)$ and the identity.
Theorem \ref{Lm:Endhomhom} then shows that there exists a homotopy
homomorphism of $\NN$-operads $\End_{S_*(D_2)}\longrightarrow
\End_{S_*(F)}$. By composition with the morphism
$\mathrm{PsOpd}\longrightarrow \End_{S_*(D_2)}$, 
the $\NN$-collection $S_*(F)$ is a homotopy algebra for the
$\NN$-operad $\mathrm{PsOpd}$ (cf. Proposition
\ref{Tm:Linftymorphofdeformcplx}). 
Both $i$ and $r$ (and thus $S_*(i)$ and $S_*(r)$) are 
compatible with the symmetric group actions on $D_2(n)$ and $F(n)$. 
Consequently, this makes the singular chains $S_*(F)$ an equivariant 
homotopy $\mathrm{PsOpd}$-algebra, and thus an operad up to
homotopy.
\end{Pf}

\maketitle\bibliographystyle{plain}
\bibliography{hopf}

~\\\textsc{Pepijn van der Laan} (\texttt{pvanderlaan@crm.es} or \texttt{vdlaan@math.uu.nl})\\
Centre de Recerca Matem\`atica, \\
Apartat 50, E-08139 Bellaterra, Spain

\end{document}